\documentclass[a4paper,10pt]{article}
\usepackage[english]{babel}
\usepackage[utf8]{inputenc}

\usepackage{amsmath}
\usepackage{amsfonts}
\usepackage{amssymb}
\usepackage{amsthm}
\usepackage{graphicx}
\usepackage{pdfpages}

\newcommand{\mytext}[1]{ \: \textrm{#1} \: }
\newcommand{\mysetdescr}[2]{\left\{ #1 \: \left| \: #2 \right. \right\} }
\newcommand{\darr}{{\downarrow \,}}
\newcommand{\uarr}{{\uparrow \,}}
\newcommand{\ouarr}{{\uparrow_{_{_{\!\!\!\!\circ}}}}}
\newcommand{\odarr}{{\downarrow^{^{\!\!\!\!\circ}}}}

\newcommand{\setx}[1]{ \{ #1 \} }

\newcommand{\myNkz}[1]{\setx{ 0, \ldots , #1 }}

\newcommand\myparal{{\|}}

\newcommand{\mf}[1]{\mathfrak{ #1 }}
\newcommand{\fL}{\mf{L}}
\newcommand{\fP}{\mf{P}}
\newcommand{\fU}{\mf{U}}

\newcommand{\LP}{{\fL(P)}}
\newcommand{\PX}{{\fP(X)}}
\newcommand{\UP}{{\fU(P)}}

\newcommand{\subseteqA}{\, \subseteq \,}

\newcommand{\cupA}{\cup \,}
\newcommand{\inA}{\in \,}

\newcommand{\Pab}{{P \setminus (a,b)}}

\newcommand{\cov}{\lessdot}

\def\BP{\begin{proof}}
\def\EP{\end{proof}}

\DeclareMathOperator{\Crit}{Crit}

\begin{document}

\theoremstyle{plain}
\newtheorem{theorem}{Theorem}
\newtheorem{definition}{Definition}
\newtheorem{corollary}{Corollary}
\newtheorem{lemma}{Lemma}
\newtheorem{proposition}{Proposition}

\title{\bf About posets for which no lower cover or no upper cover has the fixed point property}
\author{\sc Frank a Campo}
\date{\small Viersen, Germany\\
{\sf acampo.frank@gmail.com}}

\maketitle

\begin{abstract}
\noindent
For a finite non-empty set $X$, let $\PX$ denote the set of all posets with carrier $X$, ordered by inclusion of their partial order relations. We investigate properties of posets $P \in \PX$ for which no lower cover or no upper cover in $\PX$ has the fixed point property. We derive two conditions, one of them sufficient for that no lower cover of $P$ has the fixed point property, the other one sufficient for that no upper cover of $P$ has the fixed point property. If $P$ itself has the fixed point property, the conditions are even equivalent to the respective total lack of lower or upper covers with the fixed point property. We use the results to confirm a conjecture of Schr\"{o}der.
\newline

\end{abstract}

\section{Introduction}

For a finite non-empty set $X$, let $\PX$ be the set of posets with carrier $X$. We define a partial order relation $\sqsubseteq$ on $\PX$ by setting $(X,\leq) \sqsubseteq (X,\leq')$ iff $\leq \subseteqA \leq'$. Zaguia \cite{Zaguia_2008} and Schr\"{o}der \cite{Schroeder_2021} investigated the following question: Given a poset $P \in \PX$, what can we say about the fixed point property of the lower and upper covers of $P$ in $\PX$? This question is also the subject of the present article.

Zaguia \cite{Zaguia_2008} works with upper covers. He shows that a finite poset $P$ not being a chain which is dismantlable by retractables (irreducibles) has an upper cover which is dismantlable by retractables (irreducibles), too. He also proves results and provides examples for three combinations of ``to have and not to have the fixed point property'': $P$ has it and all of its upper covers have it, $P$ has it but none of its upper covers has it, $P$ does not have it but all of its upper covers have it.

Schr\"{o}der \cite{Schroeder_2021} extends the approach to lower covers. In particular, he provides examples for posets having the fixed point property for which none of the upper covers and none of the lower covers has the fixed point property. He constructs longest possible chains of posets in $\PX$ alternating in having and not having the fixed point property. Additionally, he shows that there are non-dismantlable posets $P$ having the fixed point property for which all upper covers have the fixed point property, too.

Our approach is based on \cite{Schroeder_2021}. For $P \in \PX$, let $\prec$ be the set of covering relations $a \cov b$ in $P$ in which $a$ is a minimal point of $P$ and $b$ is a maximal one. As pointed out in \cite[Theorem 1]{Schroeder_2021} (see Lemma \ref{lemma_Schroeder} below), a poset $P$ having the fixed point property becomes disconnected if we remove any edge $(a,b) \inA \prec$. Based on this observation, we say that a poset $P \in \PX$ has an {\em FPP-graph} iff it is connected but becomes disconnected if any edge $(a,b) \inA \prec$ is removed. Having an FPP-graph is thus a necessary condition for having the fixed point property, and in the case of ${\prec} = \emptyset$, the poset $P$ has an FPP-graph iff it is connected.

We call $P \in \PX$ {\em $\fL$-shielded ($\fU$-shielded)} iff no lower cover (upper cover) of $P$ has an FPP-graph. (In consequence, no lower cover (upper cover) of $P$ has the fixed point property.) Our main results are proven in Section \ref{sec_ULshielded}. In the Propositions \ref{prop_Lshielded} and \ref{prop_Ushielded}, we derive conditions for that a lower cover and an upper cover of $P$ does not have an FPP-graph. In consequence, $\fL$-shieldedness of a poset turns out to rely on a demanding property:

\begin{theorem} \label{theo_Lshielded}
If $P \in \PX$ has the {\em 3-chain-property} defined in Definition \ref{def_dreiChProp}, then $P$ is $\fL$-shielded, and if $P$ has an FPP-graph (thus, in particular, if $P$ has the fixed point property), also the inverse implication holds: $\fL$-shieldedness of $P$ implies the 3-chain-property.
\end{theorem}

The result about $\fU$-shieldedness is more straight forward, in particular in its first part. It refers to the set $\Crit(P)$ of the {\em critical pairs} of $P$:

\begin{theorem} \label{theo_Ushielded}
If $P \in \PX$ is connected with $\Crit(P) \subseteq \min P \times \max P$, then $P$ is $\fU$-shielded, and if $P$ has an FPP-graph (thus, in particular, if $P$ has the fixed point property), also the inverse implication holds: If $P$ is $\fU$-shielded, then $\Crit(P) \subseteq \min P \times \max P$.
\end{theorem}

Disconnected $\fU$-shielded posets are characterized in Proposition \ref{prop_Ushielded_disconnect}.

In Lemma \ref{lemma_BrueckenSumme} in Section \ref{sec_constructions} we show that under rather mild conditions, $\fL$-shieldedness and $\fU$-shieldedness of two disjoint posets $P$ and $Q$ is inherited by a poset which is constructed from $P$ and $Q$ in a simple way. Using this result, we are able to present a simple proof of a conjecture of Schr\"{o}der  \cite[p.\ 258]{Schroeder_2021}: for all $n \geq 37$, there exists a poset $P$ with $\# P = n$ which has the fixed point property and for which no lower cover and no upper cover has the fixed point property. However, before we found our proof, Schr\"{o}der \cite{Schroeder_Mail2024} sent us a nearby finished proof of the conjecture combining results of the author \cite{aCampo_toappear_fourCrowns} about crowns as retracts with an advanced development of the ``adding a point''-technique already used in \cite{Schroeder_2021}. Nevertheless, our proof does not rely in any way on \cite{Schroeder_Mail2024}; with the exception of using \cite{Schroeder_2021} as a common starting point, the two proofs are independent in every respect.

\section{Preparation} \label{sec_prepara5tion}

\subsection{Notation} \label{subsec_notation}

In this section, we introduce our notation and recapitulate definitions of structures being in the focus of our investigation. For all other terms of order theory, the reader is referred to standard textbooks as \cite{Schroeder_2016}. 

Let $P = (X, \leq_P)$ be a finite poset. (We deal with finite posets only in this paper.) For $Y \subseteq X$, the {\em induced sub-poset} $P \vert_Y$ of $P$ is $\left( Y, {\leq_P} \cap (Y \times Y) \right)$. To simplify notation, we identify a subset $Y \subseteq X$ with the poset $P \vert_Y$ induced by it. We call a chain consisting of $n$ points an {\em $n$-chain} and say that it has {\em length} $n-1$. The maximal length of a chain contained in $P$ is called the {\em height of $P$}. 

Let $x <_P y$. We say that $x \in X$ is a {\em lower cover} of $y \in X$ in $P$ (and $y$ an {\em upper cover} of $x$) iff $x \leq_P z \leq_P y$ implies $z \in \setx{x,y}$ for all $z \in X$. We define
\begin{align*}
\lessdot_P & := \mysetdescr{ (x,y) \inA <_P }{ x \mytext{ is a lower cover of } y }.
\end{align*}
Additionally, we define
\begin{align*}
\myparal_P & := \mysetdescr{ (x,y) \in X \times X }{ (x,y) \notin {\leq_P} \mytext{ and } (y,x) \notin {\leq_P} }.
\end{align*}

For $y \in X$, we define the {\em (punctured) down-set} and {\em (punctured) up-set induced by $y$} as
\begin{align*}
\darr_P y & := \mysetdescr{ x \in X }{ x \leq_P y }, \quad \odarr_P y := ( \darr_P y ) \setminus \setx{y}, 
\\
\uarr_P y & := \mysetdescr{ x \in X }{ y \leq_P x }, \quad \ouarr_P y := ( \uarr_P y ) \setminus \setx{y},
\end{align*}
and for $x, y \in X$, the {\em interval $[x,y]_P$} is defined as
\begin{align*}
[x,y]_P & := ( \uarr_P x ) \cap ( \darr_P y ).
\end{align*}

For $P \in \PX$, we define the following sets:
\begin{align*}
L(P) & := \mytext{the set of minimal points of} P, \\
U(P) & := \mytext{the set of maximal points of} P, \\
M(P) & := X \setminus (L(P) \cup U(P)),  \\
\mytext{and} \quad \quad
\prec_P & := \mysetdescr{ (a,b) \in {<_P} }{ M(P) \cap [ a, b ]_P = \emptyset }.
\end{align*}
For $(a,b) \inA <_P$, we have $\setx{ a, b } \subseteq [a,b]_P$, and we conclude that $\prec_P$ is a subset of $L(P) \times U(P)$. In fact, $\prec_P \; = \cov_P \cap ( L(P) \times U(P) )$, and therefore, the elements of $\prec_P$ are called {\em minmax covering relations} in \cite{Schroeder_2021}.

For disjoint posets $P = (X, \leq_P)$ and $Q = (Y, \leq_Q)$, their {\em direct sum} $P + Q$ is the poset $(X \cup Y, {\leq_P} \cup {\leq_Q})$. We get the Hasse-diagram of $P + Q$ by juxtaposing the Hasse-diagrams of $P$ and $Q$.

From the rich theory of homomorphisms and retractions, we need basic terms only. For a poset $P$ with carrier $X$, a mapping $f : X \rightarrow X$ is called an {\em endomorphism} of $P$ iff $x \leq_P y$ implies $f(x) \leq_P f(y)$ for all $x, y \in X$. A point $x \in X$ is called a {\em fixed point} of an endomorphism $f$ of $P$ iff $f(x) = x$, and we say that $P$ has the {\em fixed point property} iff every endomorphism of $P$ has a fixed point.

An endomorphism $r$ of $P$ is called a {\em retraction of $P$} iff $r$ is idempotent, and an induced sub-poset $R := P \vert_Y$, $Y \subseteq X$, is called a {\em retract of $P$} iff a retraction $r : P \rightarrow P$ exists with $r[X] = Y$.

A point $x \in P$ is called {\em irreducible} iff $\odarr_P x$ has a maximum $z$ or $\ouarr_P x$ has a minimum $z$. In both cases, the mapping $r : X \rightarrow X$ defined by 
\begin{align*}
y & \mapsto
\begin{cases}
y, & \mytext{if } y \not= x, \\
z, & \mytext{if } y = x,
\end{cases}
\end{align*}
is a retraction with retract $P \vert_{X \setminus \setx{x}}$, and we say that {\em $x$ is I-retractable to $z$}.

\subsection{Upper and lower covers of a poset} \label{subsec_poset_covers}

For the rest of this article, $X$ is a fixed finite non-empty set and the symbol $\PX$ denotes the set of all posets with carrier $X$. As already announced in the introduction, we define a partial order relation $\sqsubseteq$ on $\PX$ by setting $P \sqsubseteq Q$ iff $\leq_P \subseteqA \leq_Q$. For every $P \in \PX$, we define
\begin{align*}
\UP & := \mysetdescr{ Q \in \PX }{ Q \mytext{ is an upper cover of } P \mytext{ with respect to} \sqsubseteq }, \\
\LP & := \mysetdescr{ Q \in \PX }{ Q \mytext{ is a lower cover of } P \mytext{ with respect to} \sqsubseteq }.
\end{align*}

For $P \in \PX$, we can construct every $Q \in \PX$ with $P \sqsubseteq Q$ by successively adding pairs from $\myparal_P$ to $\leq_P$, and we can construct every $Q \in \PX$ with $Q \sqsubseteq P$ by successively removing edges from $<_P$. Therefore, we define for every $P \in \PX$ and all $(a,b) \in X \times X$ the directed graphs
\begin{align*}
P(a,b) & := \left( X, \leq_P \cup \; \setx{(a,b)} \right), \\
\Pab & := \left( X, \leq_P \setminus \setx{(a,b)} \right).
\end{align*}

According to a result of Dean and Keller \cite{Dean_Keller_1968} from 1968, the upper covers of $P \in \PX$ are characterized as follows:
\begin{align} \nonumber
\UP = & \mysetdescr{ P(a,b) }{ (a,b) \in \Crit(P) }, \\
\label{def_CritP}
\mytext{where} \quad \Crit(P) := & \mysetdescr{(a,b) \in \myparal_P }{ \odarr_P  a \subseteq \odarr_P b \; \mytext{and} \; 
\ouarr_P  b \subseteq \ouarr_P a }.
\end{align}
The elements of $\Crit(P)$ are called {\em critical pairs}. Furthermore, it is easily seen that
\begin{align*}
\LP & = \mysetdescr{ \Pab }{ (a,b) \in \cov_P }.
\end{align*}
In the following lemma, elementary facts about $P$ and $Q \in \UP$ are collected:
\begin{lemma} \label{lemma_elementFacts}
Let $P \in \PX$, $(a,b) \in \Crit(P)$, $Q := P(a,b)$. Then, for all $x \in X$,
\begin{align} \label{darr_Q_x}
\darr_Q x & =
\begin{cases}
\darr_P x \cup \setx{a}, & \mytext{if } b \in \darr_P x, \\
\darr_P x, & \mytext{otherwise};
\end{cases} \\ \label{uarr_Q_x}
\uarr_Q x & =
\begin{cases}
\uarr_P x \cup \setx{b}, & \mytext{if } a \in \uarr_P x, \\
\uarr_P x, & \mytext{otherwise};
\end{cases} \\ \label{interval_xy}
\mytext{thus} \quad 
[x,y]_P & \subseteq [x,y]_Q \subseteq [x,y]_P \cup \setx{a,b} \quad \mytext{for all } x, y \in X,
\end{align}
and $[x,y]_P = [x,y]_Q$ if $b \leq_P x$ or $y \leq_P a$. In consequence,
\begin{align}  \label{LQ_LP}
L(Q) & = L(P) \setminus \setx{b}, \\  \label{UQ_UP}
U(Q) & = U(P) \setminus \setx{a}, \\ \label{MQ_MP}
M(Q) \setminus \setx{a,b} & \subseteq M(P) \subseteq M(Q), \\
\label{precQ_precP}
\prec_Q \setminus \setx{(a,b)} & \subseteq \; \prec_P,
\end{align}
and
\begin{align} \label{precP_precQ}
\begin{split}
\prec_Q & = \; \prec_P \cup \; \setx{(a,b)} \quad \quad \quad \quad \quad \; \mytext{if } (a,b) \in L(P) \times U(P), \\ 
\prec_P \setminus \prec_Q & = \; \lessdot_P \cap \left( \min \darr_P a \times \setx{b} \right) \; \; \; \; \mytext{if } (a,b) \in M(P) \times U(P), \\ 
\prec_P \setminus \prec_Q & = \; \lessdot_P \cap \left( \setx{a} \times \max \uarr_P b \right) \quad \mytext{if } (a,b) \in L(P) \times M(P), \\
\prec_Q & = \; \prec_P \quad \quad \quad \quad \quad \quad \quad \quad \quad \; \; \mytext{in all other cases.}
\end{split}
\end{align}
\end{lemma}
\BP The first two equations are a direct consequence of the definition of $\Crit(P)$ and ${\leq_Q} = {\leq_P} \cupA \setx{(a,b)}$. The first inclusion in \eqref{interval_xy} is due to $\leq_P \subseteqA \leq_Q$; for the second one, observe that with \eqref{darr_Q_x} and \eqref{uarr_Q_x}
\begin{align*}
[x,y]_Q & \; \; \subseteq \; \; 
\left( \setx{b} \cup \uarr_P x \right) \; \cap \; \left( \setx{a} \cup \darr_P y \right) \; \; \subseteq \; \; \setx{a,b} \cup [x,y]_P.
\end{align*}
For the proof of the addendum to \eqref{interval_xy} assume $b \leq_P x$. (The proof for $y \leq_P a$ is dual.) Then $x \not\leq_P a$, hence
$$
[x,y]_Q \; = \; \uarr_Q x \cap \darr_Q y \; \stackrel{(\ref{darr_Q_x}, \ref{uarr_Q_x})}{\subseteq} \; \uarr_P x \cap ( \darr_P y \cup \setx{a}) \; = \; [x,y]_P \; \stackrel{\eqref{interval_xy}}{\subseteq} \; [x,y]_Q.
$$

The elements of $L(Q)$ are characterized by $\darr_Q x = \setx{x}$, and \eqref{LQ_LP} follows with \eqref{darr_Q_x}; \eqref{UQ_UP} is dual. In consequence, $L(Q) \cup U(Q) \subseteq L(P) \cup U(P)$, and the second inclusion in \eqref{MQ_MP} follows. Now the first one is due to the fact that according to \eqref{LQ_LP} and \eqref{UQ_UP}, the points $a$ and $b$ are the only elements of $X$ which can switch from $L(P) \cup U(P)$ to $M(Q)$.

For $(x,y) \inA \prec_Q$, \eqref{interval_xy} and \eqref{MQ_MP} yield $[x,y]_P \cap M(P) \subseteq [x,y]_Q \cap M(Q) = \emptyset$, and \eqref{precQ_precP} is shown. Now let $(x,y) \in {\prec_P} \setminus \prec_Q$. There exists an $m \in X \setminus \setx{x,y}$ with $x <_Q m <_Q y$, and because of \eqref{interval_xy}, we have $m \in \setx{a,b} \cup [x,y]_P$, hence $m \in \setx{a,b}$. For $m = a$, we have $x \in \darr_P a$, and we must have $a \in \left( \darr_Q y \right) \setminus \left( \darr_P y \right)$ which is according to \eqref{darr_Q_x} equivalent to $y \in \left( \uarr_P b \right) \setminus \left( \uarr_P a \right) \subseteq \setx{b}$ (the last inclusion is due to $(a,b) \in \Crit(P)$). The case $m = b$ is dual and all together we get
$$
\prec_P \setminus \prec_Q \; \subseteq \; \left( \min \darr_P a \times \setx{b} \right) \cup \left( \setx{a} \times \max \uarr_P b \right).
$$
The intersection of the set on the right with ${\prec_P} = ( L(P) \times U(P) ) \cap \lessdot_P$ yields the second and third equation in \eqref{precP_precQ} because the sets on the right of these equations clearly belong to ${\prec_P} \setminus {\prec_Q}$. Additionally, we get ${\prec_P} \setminus {\prec_Q} = \emptyset$ in all other cases (for $(a,b) \in L(P) \times U(P)$, observe $(a,b) \in \myparal_P$). For these cases, we have due to \eqref{precQ_precP} only to decide whether $(a,b)$ belongs to $\prec_Q$ or not:
\begin{itemize}
\item $(a,b) \in L(P) \times U(P)$: We have $a \cov_Q b$ due to
\begin{align*}
\setx{a,b} & \; \subseteq \; [a,b]_Q
\; \stackrel{\eqref{interval_xy}}{\subseteq} \; 
[a,b]_P \cup \setx{a,b} \; = \; \setx{a,b},
\end{align*}
and the equations \eqref{LQ_LP} and \eqref{UQ_UP} deliver $(a,b) \in L(Q) \times U(Q)$.
\item The points $a$ and $b$ both belong to one of the sets $L(P)$, $M(P)$, or $U(P)$: \eqref{LQ_LP}, \eqref{MQ_MP}, and \eqref{UQ_UP} show $(a,b) \notin L(Q) \times U(Q)$.
\end{itemize}

\EP

\section{$\fL$-shielded and $\fU$-shielded posets} \label{sec_ULshielded}

\begin{figure}
\begin{center}
\includegraphics[trim = 70 688 415 70, clip]{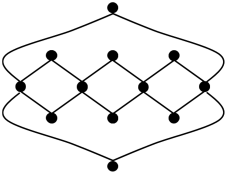}
\caption{\label{fig_BeispielLU} An example of an $\fL$- and $\fU$-shielded poset with an FPP-graph.}
\end{center}
\end{figure}

Schr\"{o}der has proven the following result:

\begin{lemma}[{\cite[Theorem 1]{Schroeder_2021}}] \label{lemma_Schroeder}
Let $P$ be a connected poset and $(a,b) \in {\prec_P}$. The poset $P$ has the fixed point property iff $\Pab$ is disconnected with both connectivity components having the fixed point property.
\end{lemma}

We pick up the constellation of connectivities and give it an own name:

\begin{definition} \label{def_FPPgraph}
We say that a poset $P$ {\em has an FPP-graph} iff $P$ is connected and $\Pab$ is disconnected for all $(a,b) \in {\prec_P}$. Furthermore, we call $P$
\begin{itemize}
\item {\em $\fL$-shielded}, iff no poset in $\LP$ has an FPP-graph;
\item {\em $\fU$-shielded}, iff no poset in $\UP$ has an FPP-graph.
\end{itemize}

\end{definition}

A simple example is the poset $P$ shown in Figure \ref{fig_BeispielLU}. It has an FPP-graph because it is connected with $\prec_P = \emptyset$, and it is $\fL$-shielded and $\fU$-shielded which can easily be seen by applying Lemma \ref{lemma_hPzwei} and Corollary \ref{coro_Ushielded} proven below.

Having an FPP-graph is a neccessary condition for having the fixed point property. A poset $P$ which is $\fL$-shielded and $\fU$-shielded has thus no lower or upper cover with the fixed point property.

A poset $P$ with ${\prec_P} = \emptyset$ has an FPP-graph iff $P$ is connected. A singleton and every disconnected poset is trivially $\fL$-shielded, and a chain $P$ is always $\fU$-shielded due to $\UP = \emptyset$. In the following sections, we derive sufficient conditions for $P$ being $\fL$-shielded or $\fU$-shielded, respectively, which are also necessary conditions if $P$ has an FPP-graph.

For a poset, being $\fL$-shielded {\em and} $\fU$-shielded is in most cases equivalent to being {\em minmax-encircled} in the sense of \cite[Definition 7]{Schroeder_2021}.

\subsection{$\fL$-shielded posets} \label{subsec_Lshielded}

\begin{figure}
\begin{center}
\includegraphics[trim = 85 665 210 70, clip]{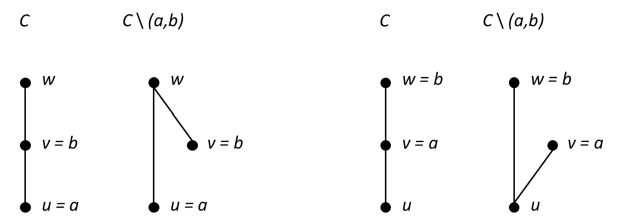}
\caption{\label{fig_Prop_1} Illustration of Definition \ref{def_dreiChProp}. There are two possible locations for the points $a$ and $b$ in the 3-chain $C$, resulting in different posets $C \setminus (a,b)$. In the first case, we have $(x,y) \in \setx{ (a,w), (b,w) }$, in the second one $(x,y) \in \setx{ (u,a), (u,b) }$. Condition \eqref{xy_LU} deals with $(x,y) \in \setx{ (a,w), (u,b) }$, whereas the conditions \eqref{x_M} and \eqref{y_M} deal with $(x,y) = (b,w)$ and $(x,y) = (u,a)$, respectively. }
\end{center}
\end{figure}

As already announced in the introduction, $\fL$-shieldedness of a poset is a demanding property. The complicated case is $P$ having an FPP-graph and $\Pab$ being connected: For $\Pab$ not having an FPP-graph, the removal of $(a,b)$ must make an edge $(x,y)$ with $(\Pab) \setminus (x,y)$ connected switching from ${<_P} \setminus \prec_P$ to $\prec_{\Pab}$. We give the required constellation an own name (for illustration, see Figure \ref{fig_Prop_1}):

\begin{definition} \label{def_dreiChProp}
Let $P \in \PX$ and $(a,b) \in \cov_P$ with $Q := \Pab$ connected. We say that $(a,b)$ has the {\em 3-chain-property} iff there exists a maximal chain $C = u \lessdot_P v \lessdot_P w$ in $P$ with $a, b \in C$ and a cover-pair $(x,y) \in \lessdot_{C \setminus (a,b)}$  for which $Q \setminus (x,y)$ is connected with 
\begin{align} \label{xy_LU}
(x,y) \in L(P) \times U(P) & \quad \Rightarrow \quad [u,w]_P = \setx{u,v,w}; \\ \label{x_M}
x \in M(P) & \quad \Rightarrow \quad b \mytext{ is I-retractable to } a \mytext{ in} P; \\ \label{y_M}
y \in M(P) & \quad \Rightarrow \quad a \mytext{ is I-retractable to } b  \mytext{ in} P.
\end{align}
Correspondingly, we say that $P$ itself has the {\em 3-chain-property} iff every cover-relation $(a,b) \in \cov_P$ with $\Pab$ connected has the 3-chain-property.
\end{definition}

Firstly, we note that Theorem \ref{theo_Lshielded} holds for every disconnected poset $P$, because such a poset has the 3-chain-property (there is no $(a,b) \in \cov_P$ with $\Pab$ connected), it is trivially $\fL$-shielded, and it does not have an FPP-graph. By treating the case ``$P$ connected'', the following proposition completes the proof of the theorem:

\begin{proposition} \label{prop_Lshielded}
Let $P \in \PX$ be connected and $(a,b) \in \cov_P$ with $\Pab$ connected. If $(a,b)$ has the 3-chain-property, then $\Pab$ does not have an FPP-graph, and if $P$ has an FPP-graph, also the inverse implication holds: If $\Pab$ does not have an FPP-graph, then $(a,b)$ has the 3-chain-property.
\end{proposition}
\BP Assume that $(a,b)$ has the 3-chain-property, $Q := \Pab$. There exists a maximal chain $C = u \lessdot_P v \lessdot_P w$ with $a, b \in C$ and there exists a pair $(x,y) \in \lessdot_{C \setminus (a,b)}$ for which $Q \setminus (x,y)$ is connected. We show $(x,y) \inA \prec_Q$. Because of $(a,b) \inA \lessdot_P$, two cases are possible:
\begin{itemize}
\item $(a,b) = (u,v)$ (left part of Figure \ref{fig_Prop_1}): Then $w = y \in U(P) \subseteq U(Q)$ according to \eqref{UQ_UP}. For $x = a \in L(P)$, alternative \eqref{xy_LU} yields $[a,y]_P = \setx{a,b,y}$, hence $[x,y]_Q = [a,y]_Q = \setx{a,y} = \setx{x,y}$, and $(x,y) \inA \prec_Q$ follows because of $x \in L(P) \subseteq L(Q)$ according to \eqref{LQ_LP}. And for $x = b$, alternative \eqref{x_M} delivers $\odarr_P x = \odarr_P b = \setx{a}$, hence $x \in L(Q)$, and $(x,y) \in {\prec_Q}$ follows again.
\item $(a,b) = (v,w)$: Dual to the previous case, resulting in $(x,y) \inA \prec_Q$, too.
\end{itemize}

Now assume that $P$ has an FPP-graph but $Q := \Pab$ not for an $(a,b) \in \cov_P$ with $Q$ connected. There exists a pair $(x,y) \inA \prec_Q$ for which $Q \setminus (x,y)$ is connected. We have $(x,y) \inA <_P$ and $x \in L(Q)$, $y \in U(Q)$. Three cases are possible:
\begin{itemize}
\item $(x,y) \in L(P) \times U(P)$: The digraph $P \setminus (x,y)$ is connected because $Q \setminus (x,y)$ is connected, hence $(x,y) \notin \; \prec_P$. There exists thus an $m \in X$ with $x <_P m <_P y$. Because of $m \notin [x,y]_Q$, we must have  $(x,m) = (a,b)$ or $(m,y) = (a,b)$. 
\begin{itemize}
\item $(x,m) = (a,b)$: $[x,y]_Q = \setx{x,y}$ delivers $[x,y]_P = \setx{ a, b , y }$. In particular, $b \lessdot_P y$. \eqref{LQ_LP} and \eqref{UQ_UP} yield $L(P) = L(Q) \setminus \setx{b} \ni a$ and $U(P) = U(Q) \setminus \setx{a} \ni y$, and $a \lessdot_P b \lessdot_P y$ is a maximal chain in $P$.
\item $(m,y) = (a,b)$: Dual to the previous case, resulting in that the chain $x \lessdot_P a \lessdot_P b$ is a maximal chain in $P$ with $[x,b]_P = \setx{ x, a, b}$.
\end{itemize}
\item $x \in M(P)$: There exists a point $z \in X$ with $(z,x) \inA <_P \setminus <_Q$. We conclude $(z,x) = (a,b)$, hence $(x,y) \in \cov_P$ due to the remark after \eqref{interval_xy}, and additionally, $x \in L(Q)$ enforces $\setx{a} = \odarr_P x$. Therefore, $a \in L(P)$, and $b=x$ is I-retractable to $a$ in $P$. Because of $y \not= a$, \eqref{UQ_UP} yields $y \in U(P)$, and due to $(a,b) \in \lessdot_P$, the chain $a \lessdot_P b \lessdot_P y$ is a  maximal chain in $P$.
\item $y \in M(P)$: Dual to the previous case, resulting in that $a=y$ is I-retractable to $b$ in $P$ and that $x \lessdot_P a \lessdot_P b$ is a maximal chain in $P$.
\end{itemize}
Therefore, $(a,b)$ has the 3-chain-property.

\EP

Let $P$ be a poset with an FPP-graph. If the heigth of $P$ is greater than two, then $P$ contains a cover-pair $(a,b) \in \lessdot_P \cap (M(P) \times M(P))$. It is easily seen that $\Pab$ is connected and that the points $a$ and $b$ cannot be contained in a maximal 3-chain. In consequence, $P$ is not $\fL$-shielded. If the height of $P$ is one, $P$ cannot contain a crown, and $P$ is $\fL$-shielded because every $Q \in \fL(P)$ is disconnected. Finally, for $P$ being an antichain, $P$ is trivially $\fL$-shielded, regardless if it has an FPP-graph (i.e., is a singleton) or not. Under the posets with an FPP-graph, it are thus only those of height two for which the question ``$\fL$-shielded or not?'' is of interest.

\begin{lemma} \label{lemma_hPzwei}
Let $P$ be a poset of height two having an FPP-graph. The following condition implies the 3-chain-property of $P$: for all $(a,b) \in \cov_P$ with $Q := \Pab$ connected, we have
\begin{align*}
(a,b) \in L(P) \times M(P) \quad & \Rightarrow \quad \exists \; z \in U(P) \mytext{: } [a,z]_P = \setx{a,b,z}, \\
& \quad \quad \quad \quad \quad \quad 
\mytext{and } Q \setminus (a,z) \mytext{connected,} \\
(a,b) \in M(P) \times U(P) \quad & \Rightarrow \quad \exists \; z \in L(P) \mytext{: } [z,b]_P = \setx{z,a,b} \\
& \quad \quad \quad \quad \quad \quad 
\mytext{and } Q \setminus (z,b) \mytext{connected.}
\end{align*}
If $P$ does not have an I-retractable point, the specified condition is even equivalent to the 3-chain-property of $P$.
\end{lemma}
\BP Assume that the specified condition holds and let $(a,b) \in \cov_P$ with $Q := \Pab$ connected. $(a,b) \in L(P) \times U(P)$ is not possible because $P$ has an FPP-graph. For $(a,b) \in L(P) \times M(P)$ and $z \in U(P)$ as specified, we have $b \cov_P z$. Furthermore, with \eqref{interval_xy}, $[a,z]_Q \subseteq [a,z]_P = \setx{a,b,z}$, and $(a,b) \notin \; \leq_Q$ delivers $[a,z]_Q = \setx{a,z}$, thus $(a,z) \in \cov_Q$. With $x := a$, $y := z$, the pair $(x,y)$ fulfills alternative \eqref{xy_LU} in the 3-chain-condition. The proof for $(a,b) \in M(P) \times U(P)$ runs analogous.

Now assume that $P$ does not contain an I-retractable point and that $P$ has the 3-chain-property. For $(a,b) \in \cov_P$ with $Q := \Pab$ connected, there exists a maximal chain $C = u \lessdot_P v \lessdot_P w$ in $P$ with $a, b \in C$ and a cover-pair $(x,y) \in \lessdot_{C \setminus (a,b)}$ for which $Q \setminus (x,y)$ is connected and fulfills one of the implications \eqref{xy_LU}-\eqref{y_M}. Because $P$ does not contain I-retractable points, it must be implication \eqref{xy_LU}, hence $(x,y) \in L(P) \times U(P)$ with $[u,w]_P = \setx{u,v,w}$. $(a,b) \in \cov_P$ leaves the choices $a = u = x, b = v, y = w$ or $x = u, a = v, b = w = y$, and the proposition is proven.

\EP

\subsection{$\fU$-shielded posets} \label{subsec_Ushielded}

Theorem \ref{theo_Ushielded} is a direct consequence of the following result:

\begin{proposition} \label{prop_Ushielded}
Let $P \in \PX$ be connected and $(a,b) \in \Crit(P)$. In the case of $(a,b) \in L(P) \times U(P)$, the poset $P(a,b)$ does not have an FPP-graph, and if $P$ has an FPP-graph, also the inverse implication holds: If $P(a,b)$ does not have an FPP-graph, then $(a,b) \in L(P) \times U(P)$.
\end{proposition}
\BP Let $P$ be connected and $(a,b) \in \Crit(P) \cap (L(P) \times U(P))$. Equation \eqref{precQ_precP} yields $(a,b) \inA \prec_{P(a,b)}$ with $P = ( P(a,b) ) \setminus (a,b)$ connected.

Now assume that $P$ has an FPP-graph and that there exists a pair $(a,b) \in \Crit(P)$ for which $Q := P(a,b)$ does not have an FPP-graph.  With $P$, also $Q$ is connected; there exists thus an $(x,y) \in {\prec_Q}$ for which $Q \setminus (x,y)$ is connected. 

In the case of $(x,y) = (a,b)$, we have $(a,b) = (x,y) \in L(Q) \times U(Q) \subseteq L(P) \times U(P)$ according to \eqref{LQ_LP} and \eqref{UQ_UP}.

Now assume $(x,y) \not= (a,b)$. Then $(x,y) \in {\prec_P}$ according to \eqref{precQ_precP}, and because $P$ has an FPP-graph, the poset $R := P \setminus (x,y)$ is disconnected. Because of
\begin{align*}
Q \setminus (x,y) & \; = \; \big( P(a,b) \big) \setminus (x,y) \; = \; \big( P \setminus (x,y) \big)(a,b) \; = \; R(a,b)
\end{align*}
is a connected digraph, the points $a$ and $b$ must belong to different connectivity components $A$ and $B$ of $R$, say, $a \in A$ and $b \in B$; additionally, one of these connectivity components contains $x$, the other one $y$.

Assume $a \in M(P) \cup U(P)$. Because of $y \in U(Q)$, $a = y$ is not possible. There exists an $\ell \in L(P)$ with $\ell <_P a$, and $a \not= y$ yields $\ell <_R a$, hence $\ell \in A$. Additionally, $(a,b) \in \Crit(P)$ delivers $\ell <_P b$. In the case of $\ell <_P m \leq_P b$ for an $m \in M(P)$, $(x,y) \in L(P) \times U(P)$ yields $(\ell,m) \inA <_R$, $(m,b) \inA \leq_R$, hence $b \in A$. Contradiction!

Therefore, $(\ell, b) \inA \prec_P$ with $\ell \in A$ and $b \in B$. Because $A$ and $B$ are different connectivity components in $R = P \setminus (x,y)$, we must have $(\ell,b) = (x,y)$, but then $x <_Q a <_Q y$ in contradiction to $(x,y) \inA \prec_Q$.

We have thus $a \notin M(P) \cup U(P)$, hence $a \in L(P)$. In the same way, we see $b \in U(P)$, and the proof is finished.

\EP

\begin{corollary} \label{coro_Ushielded}
A connected poset $P \in \PX$ is $\fU$-shielded if for all $a, b \in P$ 
\begin{align*}
\emptyset \not= \odarr_P a \subseteq \odarr_P b \quad & \Rightarrow \quad a \leq_P b \\
\mytext{and} \quad \quad
\emptyset \not= \ouarr_P b \subseteq \ouarr_P a \quad & \Rightarrow \quad a \leq_P b.
\end{align*}
\end{corollary}
\BP The assumptions ensure $\Crit(P) \subseteq \min P \times \max P$ in Theorem \ref{theo_Ushielded}.

\EP

Using this result and Lemma \ref{lemma_hPzwei}, it is easy to construct $\fL$- and $\fU$-shielded posets having an FPP-graph. Figure \ref{fig_BeispielLU} provides a blueprint.

A disconnected poset consisting of three or more connectivity components is $\fU$-shielded. In the case of two connectivity components, the situation is more complicated:

\begin{proposition} \label{prop_Ushielded_disconnect}
Let $P \in \PX$ consist of exactly two connectivity components. $P$ is not $\fU$-shielded iff both components have an FPP-graph.
\end{proposition}
\BP Let $A$ and $B$ be the two connectivity components of $P$. For $(a,b) \in \Crit(P)$ with $a$ and $b$ both belonging to the same connectivity component, $P(a,b)$ is disconnected.

Let $\Crit^*(P)$ be the set of the pairs $(a,b) \in \Crit(P)$ with $a$ and $b$ belonging to different connectivity components of $P$. We have $\Crit^*(P) \not= \emptyset$, and for $(a,b) \in \Crit^*(P)$, the sets $\odarr_P a$ and $\odarr_P b$ belong to the respective connectivity component of $a$ and $b$. Now $\odarr_P a \subseteq \odarr_P b$ yields $\odarr_P a = \emptyset$, hence $a \in L(P)$. In the same way we see $b \in U(P)$; therefore, $\emptyset \not= \Crit^*(P) \subseteq L(P) \times U(P)$.

For $(a,b) \in \Crit^*(P)$, the poset $P(a,b)$ is connected, and according to \eqref{precQ_precP}, we have $\prec_{P(a,b)} \; = \; \prec_P \cup \; \setx{(a,b)}$. If $P \vert_A$ and $P \vert_B$ both have an FPP-graph, then $P \setminus (x,y)$ consists of three connectivity components for all $(x,y) \in \prec_P$, and $[ P(a,b)]\setminus (x,y) = [P \setminus (x,y) ](a,b)$ cannot be connected. In consequence, $P(a,b)$ has an FPP-graph. Otherwise, there exists a pair $(x,y) \in {\prec_P}$ for which $P \setminus (x,y)$ has the connectivity components $A$ and $B$, and $[ P(a,b)]\setminus (x,y) = [P \setminus (x,y) ](a,b)$ is connected.

\EP

\section{A useful construction} \label{sec_constructions}

Let $P$ and $Q$ be two disjoint connected posets, let $p \in L(P)$ and $q \in U(Q)$ be points, and let $R := (P + Q)(a,b)$ be the poset constructed by adding the edge $(p,q)$ to the direct sum of $P$ and $Q$. Due to \eqref{precQ_precP}, the poset $R$ has an FPP-graph iff both posets $P$ and $Q$ have an FPP-graph. Moreover, Lemma \ref{lemma_Schroeder} tells us that $R$ has the fixed point property iff $P$ and $Q$ both have the fixed point property. The corresponding results for $\fL$-shieldedness and $\fU$-shieldedness  are

\begin{lemma} \label{lemma_BrueckenSumme}
Assume that $P$ and $Q$ have an FPP-graph.

1) $R$ is $\fL$-shielded iff $P$ and $Q$ are both $\fL$-shielded.

2) Assume that neither $P$ nor $Q$ is a singleton and that $\odarr_P x \not= \setx{p}$ and $\ouarr_Q y \not= \setx{q}$ holds for all points $x \in P$, $y \in Q$. Then $R$ is $\fU$-shielded if $P$ and $Q$ are both $\fU$-shielded.
\end{lemma}
\BP 1) Because all three posets $P$, $Q$, and $R$ have an FPP-graph, we have according to Theorem \ref{theo_Lshielded} to show that $R$ has the 3-chain-property iff $P$ and $Q$ both have the 3-chain-property. For a cover-relation $(a,b) \in \lessdot_R$, the poset $R \setminus (a,b)$ is connected iff either $(a,b) \in \lessdot_P$ with $P \setminus (a,b)$ connected or $(a,b) \in \lessdot_Q$ with $Q \setminus (a,b)$ connected. Now it is easily seen that for $(a,b) \in \lessdot_P$ (and correspondingly for $(a,b) \in \lessdot_Q$) the pair $(a,b)$ has the 3-chain-property in $P$ iff it has the 3-chain-property in $R$ because
\begin{itemize}
\item the maximal 3-chains in $P$ are exactly the maximal 3-chains in $R$ containing a point of $P$;
\item for all $(x,y) \in {<_P}$, the digraph $(\Pab) \setminus (x,y)$ is connected iff the digraph $(R \setminus (a,b)) \setminus (x,y)$ is connected;
\item $L(P) = L(R) \cap P$ and $U(P) = U(R) \cap P$;
\item for every point $z \in M(P) = M(R) \cap P$, we have $\odarr_P z = \odarr_R z$ and $\ouarr_P z = \ouarr_R z$.
\end{itemize}

2) Let $a \in M(P) \cup U(P)$. Due to our assumption $\odarr_P a \not= \setx{p}$, the set $\odarr_R a$ contains a point of $P \setminus \setx{p}$. hence $\odarr_R a \not\subseteq \odarr_R b$ for all $b \in Q$ because of $\odarr_R b \subseteq Q \cup \setx{p}$. We conclude that $\Crit(R)$ does not contain a point from $(M(P) \cup U(P)) \times Q$. In the dual way we see that $\Crit(R)$ does not contain a point from $P \times (L(Q) \cup M(Q)))$, hence
\begin{equation*}
\Crit(R) \cap (P \times Q) \; \subseteq \; L(P) \times U(Q).
\end{equation*}

Let $a \in M(Q) \cup U(Q)$. Because of $Q \not= \setx{q}$, the set $\odarr_Q a \subseteq \odarr_R a$ contains a point of $L(Q)$ and can thus not be a subset of $\odarr_R b = \odarr_P b$ for any $b \in P$. Together with the dual result for the points in $L(P) \cup M(P)$ we get
\begin{equation*}
\Crit(R) \cap (Q \times P) \; \subseteq \; L(Q) \times U(P).
\end{equation*}

Now let $(a,b) \in \Crit(R) \cap (P \times P)$, thus
$$
a \myparal_P b, \; \;  \odarr_R a \subseteq \odarr_R b, \; \mytext{and } \ouarr_R b \subseteq \ouarr_R a.
$$
The first inclusion directly yields $\odarr_P a \subseteq \odarr_P b$, and in the case of $a \not= p$ and $b \not= p$, the second inclusion delivers $\ouarr_P b \subseteq \ouarr_P a$. $b = p$ cannot hold because $a \not= p$ implies $q \notin \ouarr_R a$. Finally, if $a = p$, then $b \not= p$, thus $\ouarr_P b \subseteq (\ouarr_R a) \cap P = \ouarr_P a$. We thus have $(a,b) \in \Crit(P)$.

$\Crit(R) \cap (Q \times Q) \subseteq \Crit(Q)$ is seen in the dual way. All together, $\Crit(R)$ must be a subset of
\begin{align*}
\Crit(P) \; \; \cup \; \; \Crit(Q) & \; \; \cup \; \; (L(P) \times U(Q)) \; \; \cup \; \; (L(Q) \times U(P))
\end{align*}
and the result follows with Theorem \ref{theo_Ushielded}.

\EP

$R$ being $\fU$-shielded and $P$ not being $\fU$-shielded is equivalent to $x = p$ for every $(x,y) \in \Crit(P) \setminus (L(P) \times U(P))$.

Using these results, we can prove Schr\"{o}der's conjecture \cite[p.\ 258]{Schroeder_2021} that for all $n \geq 37$ there exists a poset $P$ with $\# P = n$ which has the fixed point property and is $\fL$-shielded and $\fU$-shielded. 

We have to mention that for the posets we are going to discuss, a poset $P$ is $\fL$-shielded and $\fU$-shielded iff it is minmax-encircled in the sense of \cite[Definition 7]{Schroeder_2021}. The slight differences between the definitions do thus not matter here.

\begin{figure}
\begin{center}
\includegraphics[trim = 70 695 140 70, clip]{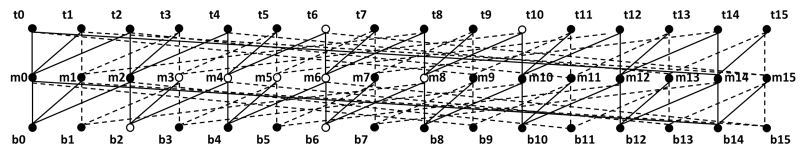}
\caption{\label{fig_X16} Schr\"{o}der's poset $X_{16}$. For better readability, the indices of point labels are not lowered and the edges are partly dotted.} 
\end{center}
\end{figure}

For $37 \leq n \leq 150$, Schr\"{o}der \cite[Theorem 2]{Schroeder_2021} proved the existence of a poset with the properties mentioned above starting with posets $X_K$ defined as follows. Let $K \geq 12$ be an even integer. The poset $X_K$ has height two and contains $3 K$ points. Its level sets are
\begin{align*}
U(X_K) & := \setx{ t_0, \ldots , t_{K-1} }, \\
M(X_K) & := \setx{ m_0, \ldots , m_{K-1} }, \\
L(X_K) & := \setx{ b_0, \ldots , b_{K-1} },
\end{align*}
and for $k \in \myNkz{K-1}$, the covering relations are
\begin{align*}
\odarr_{X_K}{m_k} & :=
\begin{cases}
\setx{ b_{k}, b_{k-2}, b_{k-7} }, & \mytext{if } k \mytext{ is even}, \\
\setx{ b_{k}, b_{k-1}, b_{k-2} }, & \mytext{if } k \mytext{ is odd},
\end{cases} \\
\ouarr_{X_K}{m_k} & :=
\begin{cases}
\setx{ t_{k}, t_{k+1}, t_{k+2} }, & \mytext{if } k \mytext{ is even}, \\
\setx{ t_{k}, t_{k+2}, t_{k+7} }, & \mytext{if } k \mytext{ is odd}.
\end{cases}
\end{align*}
(Index calculation is modulo $K$ here.) The poset $X_{16}$ is shown in Figure \ref{fig_X16}. More details about these interesting objects are contained in \cite{aCampo_toappear_fourCrowns,Schroeder_2021,Schroeder_2022_MASoC}.

Due to $\prec_{X_K} = \emptyset$, the posets $X_K$ have an FPP-graph. Furthermore, they are $\fL$-shielded and $\fU$-shielded \cite[p.\ 255--256]{Schroeder_2021}, a fact which can be checked with Lemma \ref{lemma_hPzwei} and Corollary \ref{coro_Ushielded}. Additionally they do not contain an irreducible point; they thus fulfill the assumptions made in the second part of Lemma \ref{lemma_BrueckenSumme}.

Starting with the posets $X_K$ with $12 \leq K \leq 50$, Schr\"{o}der \cite[Lemma 10]{Schroeder_2021} constructs series of super-posets of $X_K$ which are $\fL$- and $\fU$-shielded and do not contain an irreducible point. Computationally, Schr\"{o}der \cite[p.\ 258]{Schroeder_2021} checks that for every integer $n$ with $37 \leq n \leq 150$ these sequences contain a poset of cardinality $n$ which additionally has the fixed point property and in consequence also an FPP-graph.

Combining these findings with our results, it is a simple task to prove the conjecture. Let $n$ be an integer with $151 \leq n \leq 300$. From the posets with the fixed point property constructed by Schr\"{o}der, we select disjoint posets $P$ and $Q$ with $\# P + \# Q = n$ and we select points $p \in L(P)$ and $q \in U(Q)$. The poset $R := (P + Q)(p,q)$ has cardinality $n$, it has the fixed point property according to Lemma \ref{lemma_Schroeder}, and it is $\fL$-shielded and $\fU$-shielded according to Lemma \ref{lemma_BrueckenSumme}. Additionally, it does not contain an irreducible point. The conjecture is thus true for all cardinalities $n$ with $151 \leq n \leq 300$. By iterating this construction process we prove the existence of $\fL$-shielded and $\fU$-shielded posets having the fixed point property for all cardinalities $n$ with $301 \leq n \leq 600$, $601 \leq n \leq 1200$, and so on.

\end{document}